\documentclass[12pt]{article} \textwidth=6.2in \oddsidemargin=0in
\begin{document} 
\def\ov{\over} \def\t{\tau} \def\s{\sigma} \def\sp{\vspace{1ex}}
\def\cd{\cdots} \def\iy{\infty}\def\inv{^{-1}} \def\dl{\delta}
\def\be{\begin{equation}} \def\ee{\end{equation}} \def\a{\alpha}
\def\b{\beta} \def\l{\ell} \def\tn{\otimes} \def\ch{\raisebox{.3ex}{$\chi$}}
\def\ph{\varphi} \def\ps{\psi} \def\tl{\tilde} \def\pl{\partial}
\def\ga{\gamma} \def\r{\rho} \def\x{\xi} \def\eq{\equiv}
\def\ep{\varepsilon} \def\dg{{\rm diag}\,}
\def\ra{\rightarrow} \def\th{\theta} \def\arg{{\rm arg}\,}
\def\G{\Gamma} \def\D{\Delta} \def\th{\theta}

\hfill August 25, 2005

\begin{center}{\bf The Pearcey Process}\end{center}

\begin{center}{{\bf Craig A.~Tracy}\\
{\it Department of Mathematics \\
University of California\\
Davis, CA 95616, USA}}\end{center}
\begin{center}{{\bf Harold Widom}\\
{\it Department of Mathematics\\
University of California\\
Santa Cruz, CA 95064, USA}}\end{center}
\sp
\begin{abstract}
The {\it extended Airy kernel} describes the space-time correlation functions for the {\it Airy process}, which is the limiting process for a polynuclear growth model. The Airy functions themselves are given by integrals in which the exponents have a cubic singularity, arising from the coalescence of two saddle points in an asymptotic analysis. Pearcey functions are given by integrals in which the exponents have a quartic singularity, arising from the coalescence of three saddle points. A corresponding {\it Pearcey kernel} appears in a random matrix model and a Brownian motion model for a fixed time. This paper derives an {\it extended Pearcey kernel} by scaling the Brownian motion model at several times, and a system of partial differential equations whose solution determines associated distribution functions. We expect there to be a limiting nonstationary process consisting of infinitely many paths, which we call the {\it Pearcey process}, whose space-time correlation functions are expressible in terms of this extended kernel.
\end{abstract}

\renewcommand{\theequation}{1.\arabic{equation}}

\begin{center}{\bf I. Introduction}\end{center}

Determinantal processes are at the center of some recent remarkable
developments in probability theory.  These 
processes describe the mathematical structure
underpinning random matrix theory, shape fluctuations
of random Young tableaux, and certain  $1+1$ dimensional random growth
models. (See \cite{AD,FPS,JICM,soshnikov,TWICM} for recent reviews.)  Each
such process has an associated
kernel $K(x,y)$, and certain distribution functions for the process are expressed in terms of
determinants involving this kernel. (They can be ordinary determinants or operator 
determinants associated with the corresponding operator $K$ on an $L^2$ space.)
Typically these models have a parameter $n$ which might measure 
the size of the system and one is usually interested in the existence
of limiting distributions as $n\ra\iy$.  Limit laws then come down to proving that
the operator $K_n$, where we now make the $n$ 
dependence explicit, converges in trace class norm
to a limiting operator $K$.  In this context universality theorems become statements that certain 
canonical operators $K$ are the limits for a wide variety of $K_n$.  What canonical $K$ can we expect to encounter?

In various examples the kernel $K_n(x,y)$ (or, in the case of matrix kernels,
the matrix entries $K_{n,ij}(x,y)$) can be expressed as an integral
\[ \int_{C_1}\int_{C_2} f(s,t)\, e^{\phi_n(s,t;\,x,y)} \, ds\, dt .\]
To study the asymptotics of such integrals one turns to a saddle point analysis. 
Typically one finds a nontrivial limit law when there is a coalescence of saddle points.
The simplest example is the coalescence of two saddle points. This leads to the fold singularity
$\phi_2(z)={1\ov 3}\, z^3 + \lambda z$ in the theory of Thom and Arnold and a limiting kernel, the {\it Airy kernel} \cite{TWAiry}
or the more general matrix-valued 
\textit{extended Airy kernel} \cite{PS,johanssonPNG}.   

After the fold singularity comes the cusp singularity 
$\phi_3(z) ={1\ov 4}\, z^4 +\lambda_2 z^2+\lambda_1 z$.  The diffraction integrals,
which  are Airy functions
in the case of a  fold singularity, now  become Pearcey functions \cite{pearcey}.  What may be called the {\it Pearcey kernel}, since it is expressed in terms of Pearcey functions, arose in the work of Br\'ezin and Hikami \cite{BH1, BH2} on the level spacing distribution for random Hermitian matrices in an external field. More precisely, let 
$H$ be an $n\times n$ GUE  matrix (with $n$ even), suitably scaled, and $H_0$ a fixed Hermitian matrix with eigenvalues
$\pm a$ each of multiplicity $n/2$. Let $n\to\iy$. If $a$ is small the density of eigenvalues is supported in the limit 
on a single interval. If $a$ is large then it is supported on two intervals. At the ``closing of the gap'' the limiting eigenvalue distribution is described by the Pearcey kernel. 

Bleher, Kuijlaars and Aptekarev \cite{BK0,BK1,ABK2} have shown that the same kernel arises in a Brownian motion model. Okounkov and Reshetikhin \cite{RO} have encountered the same kernel in a certain growth model. 

Our starting point is with the work of Aptekarev, Bleher and Kuijlaars \cite{ABK2}.  With $n$ even again, consider $n$ nonintersecting Brownian paths
starting at position $0$ at time $\tau=0$, with half the paths conditioned
to end at $b>0$ at time $\tau=1$ and the other half conditioned to end at $-b$. 
At any fixed time this model is equivalent to the random matrix model of Br\'ezin and Hikami since they are described by the same distribution function.
If $b$ is of the order $n^{1/2}$ there is a critical time $\tau_c$
such that the limiting distribution of the Brownian paths as $n\ra\iy$ is supported by one interval for $\tau<\tau_c$ and by two intervals
when $\tau>\tau_c$. The limiting distribution at the critical time is described by the Pearcey kernel. 

It is in searching for the limiting joint distribution at several times that an {\it extended Pearcey kernel} arises.\footnote
{It was in this context that the extended Airy kernel, and other extended kernels considered in \cite{TWDyson}, arose.}
Consider times $0<\tau_1\le \cdots \le \tau_m < 1$ and ask for
the probability that for each $k$ no path passes through a set $X_k$ at
time $\tau_k$.  We show that this probability is given by the operator determinant $\det(I-K\ch)$
with an $m\times m$ matrix kernel $K(x,y)$, where
$\ch(y)=\textrm{diag}\left(\ch_{{X_k}}(y)\right)$.

We then take $b=n^{1/2}$ and scale all the times near the critical time by the substitutions $\tau_k\ra 1/2+ n^{-1/2} \t_k$ and scale the kernel by 
$x\ra n^{-1/4} x$, $y\ra n^{-1/4} y$. (Actually there are some awkward coefficients involving $2^{1/4}$ which we need not write down exactly.) The resulting limiting kernel, the extended Pearcey kernel, has
$i,j$ entry
\be -{1\ov 4\pi^2} \, \int_C\int_{-i\iy}^{i\iy} 
e^{-s^4/4+\t_j s^2/2 -y s + t^4/4 -\t_i t^2/2 + x t} \, {ds \,dt\ov s-t}\label{Pkernel}\ee
plus a Gaussian when $i<j$. The $t$ contour $C$ consists of the rays from $\pm \iy e^{i\pi/4}$ to $0$
and the rays from $0$ to $\pm e^{-i\pi/4}$.
For $m=1$ and $\t_1=0$ this reduces to the Pearcey kernel of Br\'ezin and Hikami.\footnote{In
the external source random matrix model, an interpretation is also given for the coefficients of
$s^2$ and $t^2$ in the exponential. It is not related to time as it is here.}

These authors also asked the question whether modifications of their matrix model could lead to kernels involving higher-order singularities. They found that this was so, but that the eigenvalues of the deterministic matrix $H_0$ had to be complex. Of course there are no such matrices, but the kernels describing the distribution of eigenvalues of $H_0+H$ make perfectly good sense. So in a way this was a fictitious random matrix model. In Section V we shall show how to derive analogous extended kernels and limiting processes from fictitious Brownian motion models, in which the end-points of the paths are complex numbers.

For the extended Airy kernel the authors in \cite{TWDyson} derived a system of partial differential equations, with the end-points of the intervals of the $X_k$ as independent variables, whose solution determines $\det(I-K\ch)$.\footnote
{Equations of a different kind in the case $m=2$ were found by Adler and van Moerbeke \cite{AvM}.}
Here it is assumed that each $X_k$ is a finite union of intervals.
For $m=1$ and $X_1=(\x,\iy)$ these partial differential equations reduce to ordinary differential equations which in turn can be reduced to the familiar Painlev\'e II equation.  In Section~IV of this paper we find the analogous system of partial differential equations where now the underlying kernel is the extended Pearcey kernel.\footnote
{In the case $m=1$ the kernel is integrable, i.e., it is a finite-rank kernel divided by $x-y$. (See footnote \ref{Krep1} for the exact formula.) A system of associated PDEs in this case was found in \cite{BH2}, in the spirit of \cite{TWAiry}, when $X_1$ is an interval. This method does not work when $m>1$, and our equations are completely different.} 
Unlike the case of the extended Airy kernel, here it is not until a computation at the very end that one sees that the equations close. It is fair to say that we do not really understand, from this point of view, why there should be such a system of equations.

The observant reader will have noticed that so far there has been no mention of the Pearcey process, supposedly the subject of the paper. The reason is that the existence of an actual limiting process consisting of infinitely many paths, with correlation functions and spacing distributions described by the extended Pearcey kernel, is a subtle probabilitic question which we do not now address. That for each fixed time there is a limiting random point field follows from a theorem of Lenard \cite{L1,L2} (see also \cite{soshnikov}), since that requires only a family of inequalities for the correlation functions which are preserved in the limit. But the construction of a process, a time-dependent random point field, is another matter. Of course we expect there to be one. 

\setcounter{equation}{0}\renewcommand{\theequation}{2.\arabic{equation}}

\begin{center}{\bf II. Extended kernel for the Brownian motion model}\end{center}

Suppose we have $n$ nonintersecting Brownian paths. It follows from a theorem of Karlin
and McGregor \cite{KM} that the probability density that at
times $\t_0,\ldots,\t_{m+1}$ their positions are in infinitesimal neighborhoods of
$x_{0i},\ldots,x_{m+1,i}$ is equal to
\be\prod_{k=0}^m\det\,(P(x_{m,i},x_{m+1,j},\,\s_m))\label{probdens}\ee
where 
\[\s_k=\t_{k+1}-\t_k\]
and 
\[P(x,y,\s)=(\pi\s)^{-1/2}e^{-(x-y)^2/\s}.\]

The indices $i$ and $j$ run from 0 to $n-1$, and we take
\[\t_0=0,\ \ \ \t_{m+1}=1.\]

We set all the $x_{0i}=a_i$ and $x_{m+1,j}=b_j$, thus requiring our paths to start at $a_i$ and end at $b_j$. (Later we will let all $a_i\to0$.) 

By the method of \cite{EM} (modified and somewhat simplified in \cite{TWDyson}) we shall derive an ``extended kernel'' $K(x,y)$, which is a matrix kernel
$(K_{k\l}(x,y))_{k,\l=1}^m$
such that for general functions 
$f_1,\ldots,f_m$ the expected value of
\[\prod_{k=1}^m\prod_{i=0}^{n-1}(1+f_k(x_{ki}))\]
is equal to
\[\det\,(I-K\,f),\]
where 
\[f(y)={\rm diag}\;(f_k(y)).\]
In particular the probability that for each $k$ no path passes through the set 
$X_k$ at time $\t_k$ is equal to 
\[\det\,(I-K\,\ch),\]
where
\[\ch(y)={\rm diag}\;(\ch_{X_k}(y)).\]
(The same kernel gives the correlation functions \cite{EM}. In particular the probability density (\ref{probdens}) is equal to $(n!)^{-m}\det(K_{k\l}(x_{ki},\,x_{\l j}))_{k,\l=1,\ldots,m;\;i,j=0,\ldots,n-1}$.) 

The extended kernel $K$ will be a difference $H-E$, where $E$ is the strictly upper-triangular matrix with $k,\l$ entry
$P(x,y,\,\t_\l-\t_k)$ when $k<\l$, and where $H_{k\l}(x,y)$ is given at first by the rather abstract formula (\ref{Hkl1}) below and then by the more concrete formula (\ref{Hkl2}). Then we let all $a_i\to0$ and find the integral representation (\ref{Hkl}) for the case when all the Brownian paths start at zero. This representation will enable us to take the scaling limit in the next section.

We now present the derivation of $K$.
Although in the cited references the determinants at either end (corresponding to $k=0$ and $m$ in (\ref{probdens})) were Vandermonde determinants, it is straightforward to apply the method to the present case. Therefore, rather than go through the full derivation again we just describe how one finds the extended kernel. 

For $i,\,j=0,\ldots,n-1$ we find $P_i(x)$, which are linear combination of the $P(x,a_k,\,\s_0)$ and $Q_j(y)$, which are linear combination of the $P(y,b_k,\,\s_m)$ such that
\[\int\cd\int P_i(x_1)\,\prod_{k=1}^{m-1}P(x_k,x_{k+1},\,\s_k)\,Q_j(x_m)\,dx_1\cd dx_m=\dl_{ij}.\]
Because of the semi-group property of the $P(x,y,\,\t)$ this is the same as
\be\int\int P_i(x)\,P(x,y,\,\t_m-\t_1)\,Q_j(y)\,dx\,dy=\dl_{ij}.\label{PQ}\ee

We next define for $k<\l$
\[E_{k\l}(x_k,x_\l)=\int\cdots\int \prod_{r=k}^{\l-1} P(x_r,x_{r+1},\,\s_r)\,dx_{k+1}\cdots dx_{\l-1}=P(x_k,x_\l,\;\t_\l-\t_k).\]
Set
\[P_i=P_{1i},\ \ \ \ Q_{mj}=Q_j,\]
and for $k>1$ define
\be P_{ki}(y)=\int E_{1k}(y,u)\,P_i(u)\,du=\int P(y,u,\,\t_k-\t_1)\,P_i(u)\,du,\label{Pki}\ee
and for $k<m$ define
\be Q_{kj}(x)=\int E_{km}(x,v)\,Q_j(v)\,dv=\int P(x,v,\,\t_m-\t_k)\,Q_j(v)\,dv.\label{Qkj}\ee
(These hold also for $P_{1i}$ and $Q_{mj}$ if we set $E_{kk}(x,y)=\dl(x-y)$.)

The extended kernel is given by $K=H-E$ where
\be H_{k\l}(x,y)=\sum_{i=0}^{n-1}Q_{ki}(x)\,P_{\l i}(y),\label{Hkl1}\ee
and $E_{k\l}(x,y)$ is as given above for $k<\l$ and equal to zero otherwise.

This is essentially the derivation in \cite{EM} applied to this special case. We now 
determine $H_{k\l}(x,y)$ explicitly.

Suppose 
\[ P_i(x)=\sum_k p_{ik}\,P(x,a_k,\,\s_0),\]
\[ Q_j(y)=\sum_\l q_{j\l}\,P(y,b_\l,\s_m).\]
If we substitute these into (\ref{PQ}) and use the fact that $\s_0+\t_m-\t_1+\s_m=1$ we see that it becomes
\[{1\ov\sqrt{2\pi}}\sum_{k,\l}p_{ik}\,q_{j\l}\,e^{-(a_k-b_\l)^2}=\dl_{ij}.\]
Thus, if we define matrices $P,\ Q$ and $A$ by
\[P=(p_{ij}),\ \ Q=(q_{ij}),\ \ A=(e^{-(a_i-b_j)^2}),\]
then we require $PAQ^t=\sqrt{2\pi}\,I$.

Next we compute
\[P_{ri}(y)=\int P(y,u,\,\t_r-\t_1)\,P_i(u)\,du=
\sum_k p_{ik} \,P(y,a_k,\,\t_r),\]
\[Q_{sj}(x)=\int P(x,v,\,\t_m-\t_s)\,Q_j(v)\,dv=\sum_\l q_{j\l}\,P(x,b_\l,\,1-\t_s).\]
Hence
\[H_{rs}(x,y)=\sum_i Q_{ri}(x)\,P_{si}(y)=\sum_{i,k,\l}P(x,b_\l,\,1-\t_r)\,q_{i\l}\,p_{ik}\,P(y,a_k,\,\t_s).\]
The internal sum over $i$ is equal to the $\l,k$ entry of $Q^tP=\sqrt{2\pi}\,A\inv$. So the above can be written (changing indices)
\[H_{k\l}(x,y)=\sqrt{2\pi}\,\sum_{i,j}P(x,b_i,\,1-\t_k)\,(A\inv)_{ij}\,P(y,a_j,\,\t_\l).\]
If we set $B=(e^{2\,a_i\,b_j})$ then this becomes
\be H_{k\l}(x,y)=\sqrt{2\pi}\,\sum_{i,j}P(x,b_i,\,1-\t_k)\,e^{b_i^2}\,(B\inv)_{ij}\,e^{a_j^2}\,P(y,a_j,\,\t_\l).\label{Hkl2}\ee

This gives the extended kernel when the Brownian paths start at the $a_j$. Now we are going to let all $a_j\to0$. 

There is a matrix function $D=D(a_0,\ldots,a_{n-1})$ such that for any smooth fuction $f$,
\[\lim_{a_i\to 0} D(a_0,\ldots,a_{n-1})\,\left(\begin{array}{c}f(a_0)\\f(a_1)\\\vdots\\f(a_{n-1})\end{array}\right)=\left(\begin{array}{c}f(0)\\f'(0)\\\vdots\\f^{(n-1)}(0)\end{array}\right).\]
Here $\lim_{a_i\to 0}$ is short for a certain sequence of limiting operations.
Now $B\inv$ applied to the column vector 
\[(e^{a_j^2}\,P(y,a_j,\,\t_\l))\]
equals $(DB)\inv$ applied to the vector
\[D\,(e^{a_j^2}\,P(y,a_j,\,\t_\l)).\]
When we apply $\lim_{a_i\to 0}$ this vector becomes
\[(\pl_a^j \,e^{a^2/2}\,P(y,a,\,\t_\l)|_{a=0}),\]
while $DB$ becomes the matrix $((2\,b_j)^i)$, which is invertible when all the $b_j$ are distinct. If we set $V=({b_j}^i)$ then the limiting $(DB)\inv$ is equal to $V\inv\,{\rm diag}\,(2^{-j})$.
Thus we have shown that when all $a_i=0$
\be H_{k\l}(x.y)=\sqrt{2\pi}\,\sum_{i,j}P(x,b_i,\,1-\t_k)\,e^{b_i^2}\,(V\inv)_{ij}\,2^{-j}\, \pl_a^j\,e^{a^2}\,P(y,a,\,\t_\l)|_{a=0}.\label{Hklsum}\ee

The next step is to write down an integral representation for the last factor. We have
\[e^{a^2}\,P(y,a,\,\t_\l)={1\ov \pi i\sqrt{2(1-\t_\l)}}\,e^{{y^2\ov 1-\t_\l}}\int_{-i\iy}^{i\iy} e^{{\t_\l\ov 1-\t_\l}\,s^2+2s\left(a-{y\ov1-\t_\l}\right)}\,ds.\]
Hence
\be 2^{-j}\,\pl_a^j \,e^{a^2}\,P(y,a,\,\t_\l)|_{a=0}={1\ov \pi i\sqrt{2(1-\t_\l)}}\,e^{{y^2\ov 1-\t_\l }}\int_{-i\iy}^{i\iy} s^j\,e^{{\t_\l\ov 1-\t_\l}\,s^2-{2sy\ov1-\t_\l}}\,ds.\label{Pintrep}\ee

Next we are to multiply this by $(V\inv)_{ij}$ and sum over $j$. The index $j$ appears in (\ref{Pintrep}) only in the factor $s^j$ in the integrand, so what we want to compute is
\be\sum_{j=0}^{n-1}(V\inv)_{ij}\,s^j.\label{Vsum}\ee
Cramer's rule in this case tells us that the above is equal to $\D_i/\D$,
where $\D$ denotes the Vandermonde determinant of $b=\{b_0,\ldots,b_{n-1}\}$ and $\D_i$ the Vandermonde determinant of $b$ with $b_i$ replaced by $s$. This is equal to 
\[\prod_{r\ne i}{s-b_r\ov b_i-b_r}.\]

Observe that this is the same as the residue
\be{\rm res}\left(\prod_r {s-b_r\ov t-b_r}\,{1\ov s-t},\,t=b_i\right).\label{res}\ee
This allows us to replace the sum over $i$ in (\ref{Hklsum}) by an integral over $t$. In fact, 
using (\ref{Pintrep}) and the identification of (\ref{Vsum}) with (\ref{res}) we see that the right side of (\ref{Hklsum}) is equal to
\[-{1\ov 2\pi}{1\ov \sqrt{\pi(1-\t_\l)}}e^{{y^2\ov 1-\t_\l}}\,\int_C\int_{-i\iy}^{i\iy} P(x,t,\,1-\t_k)\,e^{t^2}\,e^{{\t_\l\ov 1-\t_\l}\,s^2-{2sy\ov1-\t_\l}}\,\prod_r{s-b_r\ov t-b_r}\,{ds\,dt\ov s-t},\]
where the contour of integration $C$ surrounds all the $b_i$ and lies to one side (it doesn't matter which) of the $s$ contour. Thus 
\[H_{k\l}(x,y)=-{1\ov 2\pi^2}{1\ov \sqrt{(1-\t_k)\,(1-\t_\l)}}\,e^{{y^2\ov 1-\t_\l}-{x^2\ov 1-\t_k }}\]
\be\times\int_C\int_{-i\iy}^{i\iy} e^{-{\t_k \ov 1-\t_k}\,t^2+{2xt\ov 1-\t_k}+{\t_\l\ov 1-\t_\l}\,s^2-{2sy\ov1-\t_\l}}\;\prod_r{s-b_r\ov t-b_r}\;{ds\,dt\ov s-t}.\label{Hkl}\ee

In this representation the $s$ contour (which  passes to one side of the closed $t$ contour) may be replaced by the imaginary axis and $C$ by the contour consisting of the rays from $\pm\iy e^{i\pi/4}$  to 0 and the rays from 0 to $\pm\iy e^{-i\pi/4}$. (We temporarily call this the ``new'' contour $C$.) To see this let $C_R$ denote the new contour $C$, but with $R$ replacing 
$\iy$ and the ends joined by two vertical lines (where $t^2$ has positive real part). The $t$ contour may
be replaced by $C_R$ if the $s$~contour passes to the left of it. To show that the $s$ contour may be replaced
by the imaginary axis it is enough to show that we get 0 when the $s$~contour is the interval 
$[-iR,iR]$ plus a curve from $iR$ to $-iR$ passing around to the left of $C_R$. If we integrate first with respect to $s$ we get a
pole at $s=t$, and the resulting $t$ integral is zero because the integrand is analytic inside $C_R$.
So we can replace the $s$ contour by the imaginary axis. We then let $R\to\iy$ to see that $C_R$ may be replaced by the new contour $C$. 

\sp 

\setcounter{equation}{0}\renewcommand{\theequation}{3.\arabic{equation}}

\begin{center}{\bf III. The extended Pearcey kernel}\end{center}

The case of interest here is where half the $b_r$ equal $b$ and half equal $-b$. 
It is convenient to replace $n$ by $2n$, so that the product in the integrand in (\ref{Hkl})  is equal to 
\[\left({s^2-b^2\ov t^2-b^2}\right)^n.\]

We take the case $b=n^{1/2}$. We know from \cite{ABK2} that the critical time (the time when the support
of the limiting density changes from one interval to two) is $1/2$, and the place (where
the intervals separate) is $x=0$. We make the replacements
\[\t_k\to 1/2+n^{-1/2}\,\t_k\]
and the scaling 
\[x\to n^{-1/4}\,x,\ \ \ y\to n^{-1/4}\,y.\]
More exactly, we define
\[K_{n,ij}(x,y)=n^{-1/4}\,K_{ij}(n^{-1/4}\,x,\,n^{-1/4}\,y),\]
with the new definition of the $\t_k$. (Notice the change of indices from $k$ and $\l$ to $i$
and $j$. This is for later convenience.) 

The kernel $E_{n,ij}(x,y)$ is exactly the same as $E_{ij}(x,y)$.
As for $H_{n,ij}(x,y)$, its integral representation is obtained from (\ref{Hkl}) by the scaling replacements and then by the substitutions $s\to n^{1/4}\,s,\ t\to n^{1/4}\,t$ in the integral itself. The result is (we apologize for its length)
\[ H_{n,ij}(x,y)=-{1\ov \pi^2}{1\ov \sqrt{(1-2n^{-1/2}\t_i)\,
(1-2n^{-1/2}\t_j)}}\exp\left\{{2y^2\ov n^{1/2}-2\t_j}-
{2x^2\ov n^{1/2}-2\t_i }\right\}\]
\[\times\int_C\int_{-i\iy}^{i\iy} \exp\left\{-n^{1/2}\,{1+2n^{-1/2}\t_i \ov 1-2n^{-1/2}\t_i}\,t^2+{4xt\ov 1-2n^{-1/2}\t_i}\right\}\]
\be\times\exp\left\{n^{1/2}\,{1+2n^{-1/2}\t_j \ov 1-2n^{-1/2}\t_j}\,s^2-{4ys\ov 1-2n^{-1/2}\t_ji}\right\}\;
\left({1-s^2/n^{1/2}\ov 1-t^2/n^{1/2}}\right)^n\;{ds\,dt\ov s-t}.\label{Hnij}\ee

We shall show that this has the limiting form
\be H^{{\rm Pearcey}}_{ij}(x,y)=-{1\ov \pi^2}\int_C \int_{-i\iy}^{i\iy}e^{-s^4/2+4\,\t_j\,s^2-4ys+
t^4/2-4\,\t_i\,t^2+4xt}\,{ds\,dt\ov s-t},\label{Pkernel2}\ee
where, as in (\ref{Hkl}) and (\ref{Hnij}), the $s$ integration is along the imaginary axis and the  contour $C$ consists of the rays from $\pm\iy e^{i\pi/4}$  to 0 and the rays from 0 to $\pm\iy e^{-i\pi/4}$. 
Precisely, we shall show that
\be\lim_{n\to\iy}H_{n,ij}(x,y)=H^{{\rm Pearcey}}_{ij}(x,y)\label{limit}\ee
uniformly for $x$ and $y$ in an arbitrary bounded set, and similarly
for all partial derivatives.\footnote
{The constants in (\ref{Pkernel2}) are different from those in (\ref{Pkernel}), a matter of no importance. In the next section we shall make the appropriate change so that they agree.}

The factor outside the integral in (\ref{Hnij}) converges to $-1/\pi^2$. The first step in proving the convergence of the integral in (\ref{Hnij}) to that in (\ref{Pkernel2}) is to establish pointwise convergence of the integrand. 

The first exponential factor in the integrand in (\ref{Hnij}) is
\[\exp\Big\{-(n^{1/2}+4\t_i+O(n^{-1/2}))\,t^2+(4+O(n^{-1/2}))\,xt\Big\},\]
while the second exponential factor is
\be\exp\Big\{(n^{1/2}+4\t_j+O(n^{-1/2}))\,s^2-(4+O(n^{-1/2}))\,ys\Big\}.\label{sest}\ee
When $s=o(n^{1/4}),\ t=o(n^{1/4})$ the last factor in the integrand is equal to
\[\exp\Big\{n^{1/2}\,t^2+t^4/2+o(t^4/n)-n^{1/2}\,s^2-s^4/2+o(s^4/n)\Big\}.\]
Thus the entire integrand (aside from the factor $1/(s-t)$) is
\[\exp\Big\{-(1+o(1))\,s^4/2+(4\t_j+O(n^{-1/2}))\,s^2-(4+O(n^{-1/2}))\,ys\}\]
\be\times \exp\Big\{(1+o(1))\,t^4/4-(4\t_i+O(n^{-1/2}))\,t^2+(4+O(n^{-1/2}))\,xt\Big\}\label{test}\ee

In particular this establishes pointwise convergence of the integrands in (\ref{Hnij}) to that in 
(\ref{Pkernel2}). For the claimed uniform convergence of the integrals and their derivatives
it is enough to show that they all converge pointwise and boundedly. 
To do this we change the $t$ contour $C$ by rotating its rays slightly toward the real line. (How much we rotate we say below. We can revert to the original contour after 
taking the limit.) This is so that on the modified contour, which we denote by $C'$, we have $\Re \,t^2>0$ as well as $\Re\, t^4<0$.

The function
$1/(s-t)$ belongs to $L^q$ for any $q<2$ in the neighborhood of $s=t=0$ on the contours of integration and to $L^q$ for any $q>2$ outside this neighborhood. To establish pointwise bounded convergence of the integral it therefore suffices to show that for any $p\in(1,\iy)$ the rest of the integrand (which we know converges pointwise) has $L^p$ norm which is uniformly bounded in $x$ and $y$.\footnote
{That this suffices follows from the fact, an exercise, that if $\{f_n\}$ is a bounded sequence in $L^p$ converging pointwise
to $f$ then $(f_n,\,g)\to (f,\,g)$ for all $g\in L^q$, where $p=q/(q-1)$. We take $f_n$ to be the integrand in (\ref{Hnij}) except for the factor $1/(s-t)$ and $g$ to be $1/(s-t)$, and apply this twice, with $q<2$ in a neighborhood of $s=t=0$ and with $q>2$ outside the neighborhood.}
The rest of the integrand is the product of a function of $s$ and a function of $t$ and we show that both of these functions have uniformly bounded $L^p$ norms. 

From (\ref{test}) it follows that for some small $\ep>0$ the function of $t$ is $O(e^{\Re\,t^4/2+O(|t|^2)})$ uniformly in $x$ if $|t|\le \ep\,n^{1/4}$. Given this $\ep$ we choose $C'$ to consist of the rays of $C$ rotated slightly toward the real axis so that if $\th={\rm arg}\,t^2$ when $t\in C'$ then $\cos\th=\ep^2/2$. 

When $t\in C'$ and $|t|\le \ep\,n^{1/4}$ the function of $t$ is $O(e^{\Re\,t^4/2+O(|t|^2)})=O(e^{\cos2\th\,|t|^4/2})$. Since $\cos2\th<0$ the $L^p$ norm on this part of $C'$ is $O(1)$. 

When $t\in C'$ and $|t|\ge \ep\,n^{1/4}$ we have 
\[|1-t^2/n^{1/2}|^2=1+n\inv |t|^4-2n^{-1/2}\cos\th\,|t|^2.\]
But 
\[n\inv |t|^4-2n^{-1/2}\cos\th\,|t|^2=n^{-1/2}|t|^2\,(n^{-1/2}|t|^2-\ep^2)\ge0\]
when $|t|\ge\ep n^{1/4}$. Thus $|1-t^2/n^{1/2}|\ge1$ and the function of $t$ is $O(e^{-\cos\th\,n^{1/2}\,|t|^2+O(|t|^2)})=O(e^{-\cos\th\,n^{1/2}\,|t|^2/2})$, and the $L^p$ norm on this part of $C'$ is $o(1)$.

We have shown that on $C'$ the function of $t$ has uniformly bounded $L^p$ norm.

For the $L^p$ norm of the function of $s$ we see from (\ref{sest}) that the integral of its $p$th power is at most a constant independent of $y$ times
\[\int_0^\iy e^{-pn^{1/2}s^2+\t s^2}\,(1+n^{-1/2}\,s^2)^{pn}\,ds.\]
(We replaced $s$ by $is$, used evenness, and took any $\t>-4p\,\t_j$.) The variable change 
$s\to n^{1/4}\,s$
replaces this by
\[n^{1/4}\int_0^\iy e^{-pn\,(s^2-\log\,(1+s^2))+n^{1/2}\,\t\,s^2}\,ds.\]
The integral over $(1,\iy)$ is exponentially small. Since $s^2-\log\,(1+s^2)\ge s^4/2$ when $s\le1$, what remains is at most
\[n^{1/4}\int_0^1 e^{-pn\,s^4/2+O(n^{1/2}\,s^2)}\,ds
= \int_0^{n^{1/4}} e^{-p\,s^4/2+O(s^2)}\,ds,\]
which is $O(1)$.

This completes the demonstration of the bounded pointwise convergence of $H_{n,ij}(x,y)$
to $H^{{\rm Pearcey}}_{ij}$.
Taking any partial derivative just inserts in the integrand a polynomial in $x,\, y,\, s$ and $t$, 
and the argument for the modified integral is virtually the same.  This completes the proof of (\ref{limit}).

\setcounter{equation}{0}\renewcommand{\theequation}{4.\arabic{equation}}

\begin{center}{\bf IV. Differential equations for the Pearcey process}\end{center}

We expect the extended Pearcey kernel to characterize the Pearcey process, a point process which can be thought of 
as a infinitely many nonintersecting paths. Given sets $X_k$, the probability that for each $k$
no path passes through the set $X_k$ at time $\t_k$ is equal to 
\[\det\,(I-K\,\ch),\]
where
\[\ch(y)={\rm diag}\;(\ch_{X_k}(y)).\]

The following discussion follows closely that in \cite{TWDyson}. We take the case where each $X_k$
is a finite union of intervals with end-points $\x_{kw},\ w=1,\ 2,\ldots$, in increasing order. If we set
\[R=(I-K\,\ch)\inv\,K,\]
with kernel $R(x,y)$, then
\[\pl_{kw}\det\,(I-K\,\ch)=(-1)^{w+1}R_{kk}(\x_{kw},\x_{kw}).\]
(We use the notation $\pl_{kw}$ for $\pl_{\x_{kw}}$.) We shall find a system of PDEs in the 
variables $\x_{kw}$ with the right sides above among
the unknown functions. 

In order to have the simplest coefficients later we make the further variable changes
\[s\to s/2^{1/4},\ \ \ t\to t/2^{1/4}\]
and substitutions
\[x\to 2^{1/4}x/4,\ \ \ y\to 2^{1/4}y/4,\ \ \ \t_k\to2^{1/2}\t_k/8.\]
The resulting rescaled kernels are (we omit the superscripts ``Pearcey'')
\[H_{ij}(x,y)=-{1\ov 4\pi^2}\int_C \int_{-i\iy}^{i\iy}e^{-s^4/4+\t_j\,s^2/2-ys+
t^4/4-\t_i\,t^2/2+xt}\,{ds\,dt\ov s-t},\]
which is (\ref{Pkernel}), and
\[E_{ij}(x,y)={1\ov\sqrt{2\pi\,(\t_j-\t_i)}}e^{-{(x-y)^2\ov2(\t_j-\t_i)}}.\]

Define the vector functions
\[\ph(x)=\left({1\ov 2\pi i}\int_Ce^{t^4/4-\t_k\,t^2/2+xt}dt\right),\ \ \ 
\ps(y)=\left({1\ov 2\pi i}\int_{-i\iy}^{i\iy}
e^{-s^4/4+\t_k\,s^2/2-ys}ds\right).\]
We think of $\ph$ as a
column $m$-vector and $\ps$ as a row $m$-vector. Their components are Pearcey functions.\footnote{\label{Krep1}
In case $m=1$ the kernel has the explicit representation \[K(x,y)={\ph''(x)\,\ps(y)-\ph'(x)\,\ps'(y)+\ph(x)\,\ps''(y)-\t\,\ph(x)\ps(y)\ov x-y}\]
in terms of the Pearcey functions. (Here we set $\t=\t_1$. This is the same as the $i,i$ entry of the matrix kernel if we set $\t=\t_i$.) This was shown in \cite{BH1}. Another derivation will be given in footnote \ref{Krep2}.} The vector functions satisfy the differential equations
\[\ph'''(x)-\t\,\ph'(x)+x\,\ph(x)=0,
\ \ \ \ps'''(y)-\ps'(y)\,\t-y\,\ps(y)=0,\]
where $\t={\rm diag}\;(\t_k)$.

Define the column vector function $Q$ and row vector function $P$ by
\be Q=(I-K\,\ch)\inv\ph,\ \ \ P=\ps\,(I-\ch\,K)\inv.\label{QandP}\ee
The unknowns in our equations will be six vector functions indexed by the end-points $kw$ of the
$X_k$ and three matrix functions with the same indices. The vector functions are denoted by
\[q,\ q',\ q'',\ p,\ p',\ p''.\]
The first three are column vectors and the second three are row vectors. They are defined by
\[q_{iu}=Q_i(\x_{iu}),\ \ \ q'_{iu}=Q_i'(\x_{iu}),\ \ \ q''_{iu}=Q''_i(\x_{iu}),\]
and analogously for $p,\ p',\ p''$. The matrix function unknowns are
\[r,\ r_x,\ r_y\]
defined by
\[r_{iu,jv}=R_{ij}(\x_{iu},\,\x_{jv}),\ \ \ r_{x,iu,jv}=R_{xij}(\x_{iu},\,\x_{jv}),\ \ \
r_{y,iu,jv}=R_{yij}(\x_{iu},\,\x_{jv}).\]
(Here $R_{xij}$, for example, means $\pl_x R_{ij}(x,y)$.)

The equations themselves will contain the matrix functions $r_{xx},\ r_{xy},\ r_{yy}$ 
defined analogously, but we shall see that the combinations of them that appear
can be expressed in terms of the unknown functions.

The equations will be stated in differential form. We use the notation
\[\x=\dg(\x_{kw}),\ \ \ d\x=\dg(d\x_{kw}),\ \ \ s=\dg ((-1)^{w+1}).\]
Recall that $q$ is a column vector and $p$ a row vector.

Our equations are
\begin{eqnarray}
dr&=&-r\,s\,d\x\,r+d\x\,r_x+r_y\,d\x,\label{de1}\\
dr_x&=&-r_x\,s\,d\x\,r+d\x\,r_{xx}+r_{xy}\,d\x,\label{de2}\\
dr_y&=&-r\,s\,d\x\,r_y+d\x\,r_{xy}+r_{yy}\,d\x,\label{de3}\\ 
dq&=&d\x\,q'-r\,s\,d\x\,q,\label{de4}\\
dp&=&p'\,d\x-p\,d\x\,s\,r,\label{de5}\\
dq'&=&d\x\,q''-r_x\,s\,d\x\,q,\label{de6}\\
dp'&=&p''\,d\x-p\,d\x\,s\,r_y,\label{de7}\\
dq''&=&d\x\,(\t\,q'-\x\,q+r\,s\,q''-r_y\,s\,q'+r_{yy}\,s\,q-r\,\t\,s\,q)-r_{xx}\,s\,d\x\,q,\label{de8}\\
dp''&=&(p'\,\t+p\,\x+p''\,s\,r-p'\,s\,r_x+p\,s\,r_{xx}-p\,s\,\t\,r)\,d\x-p\,d\x\,s\,r_{yy}.\label{de9}
\end{eqnarray}

One remark about the matrix $\t$ in equations (\ref{de8}) and (\ref{de9}). Earlier $\t$ was the 
$m\times m$ diagonal matrix with $k$ diagonal entry $\t_k$. In the equations here it is the 
diagonal matrix with $kw$ diagonal entry $\t_k$. The exact meaning of $\t$ when it appears
will be clear from the context.

As in \cite{TWDyson}, what makes the equations possible is that the operator $K$ has some nice
commutators. In this case we also need a miracle at the end.\footnote{That it seems a miracle to us shows that we do not really understand why the equations should close.} 

Denote by $\ph\tn\ps$ the operator with matrix kernel $(\ph_i(x)\,\ps_j(y))$, where $\ph_i$ and $\ps_j$ are
the components of $\ph$ and $\ps$, respectively.

If we apply the operator $\pl_x+\pl_y$ to the integral defining $H_{ij}(x,y)$
we obtain the commutator relation $[D,\,H]=-\ph\tn\ps$, where $D=d/dx$. Since also 
$[D,\,E]=0$ we have 
\be [D,\,K]=-\ph\tn\ps.\label{Com1}\ee
This is the first commutator. From it follows
\[[D,\,K\,\ch]=-\ph\tn\ps\,\ch+K\,\dl.\]
Here we have used the following notation: $\dl_{kw}$ is the diagonal matrix operator 
whose $k$th diagonal entry equals multiplication by $\dl(y-\x_{kw})$, and 
\[\dl=\sum_{kw}(-1)^{w+1}\dl_{kw}.\]
It appears above because $D\,\ch=\dl$.

Set
\[\r=(I-K\,\ch)\inv,\ \ \ R=\r\,K=(I-K\,\ch\,K)\inv-I.\]
It follows from the last commutator upon left- and right-multiplication by $\r$ that
\[ [D,\,\r]=-\r\,\ph\tn\ps\,\ch\,\r+R\,\dl\,\r.\]
{}From the commutators of $D$ with $K$ and $\r$ we compute
\[[D,\,R]=[D,\,\r\,K]=\r\,[D,\,K]+[D,\,\r]\,K=-\r\,\ph\tn\ps\,(I+\ch\,R)+R\,\dl\,R.\]
Notice that 
\[I+\ch\,R=(I-\ch\,K)\inv.\]
If we recall (\ref{QandP}) we see that we have shown
\be [D,\,R]=-Q\tn P+R\,\dl\,R.\label{L1}\ee 

To obtain our second commutator we observe that if we apply $\pl_t+\pl_s$ to the integrand 
in the formula for $H_{ij}(x,y)$ we get zero
for the resulting integral. If we
apply it to $(s-t)\inv$ we also get zero. Therefore we get zero if we apply it to the 
numerator, and this operation brings down the factor 
\[t^3-\t_i\,t+x-s^3+\t_j\,s-y.\]
The same factor results if we apply to $H_{ij}(x,y)$ the operator
\[\pl_x^3+\pl_y^3-(\t_i\,\pl_x+\t_j\,\pl_y)+(x-y).\]
We deduce
\[[D^3-\t D+M,\;H]=0.\]
One verifies that also $[D^3-\t D+M,\;E]=0.$ Hence
\be[D^3-\t D+M,\;K]=0.\label{Com2}\ee
This is the second commutator.\footnote{\label{Krep2}
{}From (\ref{Com1}) we obtain also $[D^3,\,K]=-\ph''\tn\ps+\ph'\tn\ps'-\ph\tn\ps''$. Combining this with (\ref{Com1}) itself and (\ref{Com2}) for $m=1$ with $\t=\t_1$ we obtain
$[M,\,K]=\ph''\tn\ps-\ph'\tn\ps'+\ph\tn\ps''-\t\,\ph\tn\ps$. This is equivalent to the formula stated in footnote \ref{Krep1}.}
From it we obtain
\[[D^3-\t D+M,\;K\,\ch]=K\,[D^3-\t D+M,\;\ch]=K\,(D\dl D+D^2\dl+\dl D^2-\t\,\dl),\]
and this gives
\be [D^3-\t D+M,\;\r]= R\,(D\dl D+D^2\dl+\dl D^2)\,\r-R\,\t\,\dl\,\r,\label{preL2}\ee
which in turn gives
\be[D^3-\t D+M,\;R]= [D^3-\t D+M,\;\r\,K]=R\,(D\dl D+D^2\dl+\dl D^2)\,R-
R\,\t\,\dl\,R.\label{L2}\ee

Of our nine equations the first seven are universal --- they do not depend on the
particulars of the kernel $K$ or vector functions $\ph$ or $\ps$. (The same was
observed in \cite{TWDyson}.) What are not universal are equations
(\ref{de8}) and (\ref{de9}). For the equations to close we shall also have to show that the
combinations of the entries of $r_{xx},\ r_{xy}$ and $r_{yy}$ which actually appear in the equations
are all expressible in terms of the unknown functions. The reader can check that these are the
diagonal entries of $r_{xx}+r_{xy}$ and $r_{xy}+r_{yy}$ (which also give the
diagonal entries of $r_{xx}-r_{yy}$) and the off-diagonal entries 
of $r_{xx},\ r_{xy}$ and $r_{yy}$. 

What we do at the beginning of our derivation will be a repetition of what was done
in \cite{TWDyson}. First, we have
\be\pl_{kw}\r=\r\,(K\,\pl_{kw}\ch)\,\r=(-1)^w\,R\,\dl_{kw}\,\r.\label{plr}\ee
{}From this we obtain $\pl_{kw}R=(-1)^w\,R\,\dl_{kw}\,R$, and so
\[\pl_{kw} \,r_{iu,jv}=(\pl_{kw}\,R_{ij})(\x_{iu},\x_{jv})+
R_{xij}(\x_{iu},\x_{jv})\,\dl_{iu,kw}+R_{yij}(\x_{iu},\x_{jv})\,\dl_{jv,kw}.\]
\[=(-1)^w\,r_{iu,kw}\,r_{kw,jv}+r_{x,iu,jv}\,\dl_{iu,kw}+r_{y,iu,jv}\,\dl_{jv,kw}.\]
Multipliying by $d\x_{kw}$ and summing over $kw$ give (\ref{de1}). 
Equations (\ref{de2}) and (\ref{de3}) are derived analogously. 

Next we derive (\ref{de4}) and (\ref{de6}). Using (\ref{plr}) applied to $\ph$ we obtain
\[\pl_{kw} \,q_{iu}=Q'_i(\x_{iu})\,\dl_{iu,kw}+(-1)^w\,(R\,\dl_{kw}\,Q)_i(\x_{iu})=
q'_{iu}\,\dl_{iu,kw}+(-1)^w\,r_{iu,kw}\,q_{kw}.\]
Multiplying by $d\x_{kw}$ and summing over $kw$ give (\ref{de4}). If we multiply
(\ref{plr}) on the left by $D$ we obtain $\pl_{kw}\r_x=(-1)^w\,R_x\,\dl_{kw}\,\r$ and 
applying the result to $\ph$ we obtain (\ref{de6}) similarly.

For (\ref{de8}) we begin as above, now applying $D^2$ to (\ref{plr}) on the left and $\ph$
on the right to obtain
\be\pl_{kw} \,q''_{iu}=Q'''_i(\x_{iu})\,\dl_{iu,kw}+(-1)^w\,r_{xx,iu,kw}\,q_{kw}.\label{plq''}\ee
Now, though, we have to compute the first term on the right. To do it we apply (\ref{preL2})
to $\ph$ and use the differential equation satisfied by $\ph$ to obtain
\be Q'''(x)-\t\,Q'(x)+x\,Q(x)=-R_y\,\dl\,Q'+R_{yy}\,\dl\,Q+R\,\dl\,Q''-R\,\t\,\dl\,Q.
\label{Q'''}\ee
This gives 
\[Q'''_i(\x_{iu})=\t_i\,q'_{iu}-\x_{iu}\,q_{iu}+
(-r_y\,s\,q'+r_{yy}\,s\,q+r\,s\,q''-r\,\t\,s\,q)_{iu}.\]
If we substitute this
into (\ref{plq''}), multiply by $d\x_{kw}$ and sum over $kw$ we obtain equation (\ref{de8}).

This completes the derivation of the equations for the differentials of $q,\ q'$ and $q''$.
We could say that the derivation of the equations for the differentials of $p,\ p'$ and $p''$
is analogous, which is true. But here is a better way. Observe that the  $P$ for the operator
$K$ is the transpose of the $Q$ for the transpose of $K$, and similarly with $P$ and $Q$
interchanged. It follows from this that for any equation involving $Q$ there is another
one for $P$ obtained by replacing $K$ by its transpose (and so interchanging $\pl_x$ and 
$\pl_y$) and taking transposes. The upshot is that equations (\ref{de5}), (\ref{de7}) and (\ref{de9})
are consequences of (\ref{de4}), (\ref{de6}) and (\ref{de8}). The reason for the difference in
signs in the appearance of $\x$ on the right sides of (\ref{de8}) and (\ref{de9}) is the difference in 
signs in the last terms in the differential equations for $\ph$ and $\ps$.

Finally we have to show that 
the diagonal entries of $r_{xx}+r_{xy}$ and $r_{xy}+r_{yy}$  
and the off-diagonal entries of $r_{xx},\ r_{xy}$ and $r_{yy}$ are all known, in the sense
that they are expressible in terms of the unknown functions. This is really the heart of
the matter.

We use $\eq$ between expressions involving $R,\ Q$ and $P$ and their derivatives to indicate that the 
difference involves at most two derivatives of $Q$ or $P$ and at most one derivative of $R$. The reason
is that if we take the appropriate entries evaluated at the appropriate points we obtain a known quantity,
i.e., one expressible in terms of the unknown functions.

If we multiply (\ref{L1}) on the left or right by $D$ we obtain
\[R_{xx}+R_{xy}=-Q'\tn P+R_x\,\dl\,R,\ \ \ R_{xy}+R_{yy}=-Q\tn P'+R\,\dl\,R_y.\]
In particular
\[R_{xx}+R_{xy}\eq0,\ \ \ R_{xy}+R_{xy}\eq0,\]
so in fact all entries of $r_{xx}+r_{xy}$ and $r_{xy}+r_{yy}$ are known.

{}From (\ref{L1}) we obtain consecutively
\be [D^2,\,R]=-Q'\tn P+Q\tn P'+DR\,\dl\,R+R\,\dl\,R D,\label{D2com}\ee
\be [D^3,\,R]=-Q''\tn P+Q'\tn P'-Q\tn P''
+D^2R\,\dl\,R+DR\,\dl R D+R\,\dl\,R D^2.\label{D3com}\ee
If we subtract (\ref{L2}) from (\ref{D3com}) we find
\[[\t\,D-M,\;R]= -Q'' \tn P +Q' \tn P' -Q \tn P'' \]
\[+R_{xx}\,\dl\,R-R_x\,\dl R_y+R\,\dl\,R_{yy}
+R_y\,\dl R_x-R_{yy}\dl\,R-R\,\dl R_{xx}+R\,\t\,\dl\,R.\]
We use
\[R_{xx}-R_{yy}=Q\tn P'-Q'\tn P+R_x\,\dl\,R-R\,\dl\,R_y\]
and 
\[R_y=-R_x-Q\tn P+R\,\dl\,R\]
to see that this equals
\pagebreak
\[-Q'' \tn P +Q' \tn P' -Q \tn P''+(Q\tn P'-Q'\tn P-R\,\dl\,R_y)\,\dl\,R\]
\be -R\,\dl\,(Q\tn P'-Q'\tn P+R_x\,\dl\,R-R\,\dl\,R_y)+R_x\,\dl\,Q\tn P\label{form1}\ee
\[-(Q\tn P-R\,\dl\,R)\,\dl\,R_x+R\,\t\,\dl\,R.\]

We first apply $D$ acting on the left to this, and deduce that
\[D\,[\t\,D,\;R]\eq -Q''' \tn P+R_{xx}\,\dl\,Q\tn P.\]
If we use (\ref{Q'''}) and the fact that $R_{yy}\eq R_{xx}$ we see that this is $\eq0$.
This means that 
\[\t_i\,R_{xx,i,j}+\t_j\,R_{xy,i,j}\eq 0.\]
Since $R_{xx,i,j}+R_{xy,i,j}\eq 0$ we deduce that $R_{xx,i,j}$ and 
$R_{xy,i,j}$ are individually $\eq 0$ when $i\ne j$. Therefore $r_{xx,iu,jv}$
and $r_{xy,iu,jv}$ are known then.

But we still have to show that $r_{iu,iv}$ is known when $u\ne v$, and for this we apply
$D^2$ to (\ref{form1}) rather than $D$. We get this time
\[D^2\,[\t\,D-M,\;R]\eq-Q'''' \tn P +Q'''\tn (P' -P\,\dl\,R)-R_{xx}\,\dl\,R_y\,\dl\,R\]
\be -R_{xx}\,\dl\,(Q\tn P'-Q'\tn P+R_x\,\dl\,R-R\,\dl\,R_y)\label{form2}\ee
\[+R_{xxx}\,\dl\,Q\tn P+R_{xx}\,\dl\,R\,\dl\,R_x+R_{xx}\,\t\,\dl\,R.\]

We first compare the diagonal entries of $D^2\,[\t\,D,\,R]$ on the left with those of 
$R_{xx}\,\t\,\dl\,R$ on the right. The diagonal entries of the former are the same as those
of
\[\t\,(R_{xxx}+R_{xxy})\eq \t\,R_{xx}\,\dl\,R.\]
(Notice that applying $D^2$ to (\ref{L1}) on the left gives $R_{xxx}+R_{xxy}\eq R_{xx}\,\dl\,R$.) 
The difference between this and $R_{xx}\,\t\,\dl\,R$ is $[\t,\,R_{xx}]\,\dl\,R$.
Only the off-diagonal entries of $R_{xx}$ occur here, so this is $\eq0$.

If we remove these terms from (\ref{form2}) the left side becomes $-D^2\,[M,\,R]$
and the resulting right side we write, after using the fact $R_x+R_y=-Q\tn P+R\,\dl\,R$
twice, as
\[-Q'''' \tn P +Q''' \tn (P'-P\,\dl\,R)\]
\be -R_{xx}\,\dl\,(Q\tn P'-Q'\tn P-Q\tn P\,\dl\,R+R\,\dl\,Q\tn P)
+R_{xxx}\,\dl\,Q\tn P.\label{form3}\ee

Now we use (\ref{Q'''}) and the facts 
\[R_{yy}\eq R_{xx},\ \ \ R_{xy}\eq -R_{xx},\ \ \ R_{xxx}-R_{xyy}\eq R_{xx}\dl\,R\] 
(the last comes from applying $D$ to (\ref{D2com}) on the left) to obtain
\[Q'''\eq R_{xx}\,\dl\,Q,\ \ \ Q''''\eq R_{xx}\,\dl\,Q'+(R_{xxx}-R_{xx}\,\dl\,R)\,\dl\,Q.\]
Substituting these expressions for $Q'''$ and $Q''''$ into (\ref{form3})
shows that it is $\eq0$. 

This was the miracle.

We have shown that
\[D^2\,[M,\,R]\eq 0,\]
in other words $(x-y)\,R_{xx}(x,\,y)\eq 0$. If we set $x=\x_{iu},\ y=\x_{iv}$ we 
deduce that $R_{xx}(\x_{iu},\x_{iv})=r_{xx,iu,iv}$ is known when $u\ne v$. 

\sp
\setcounter{equation}{0}\renewcommand{\theequation}{5.\arabic{equation}}
\begin{center}{\bf V. Higher-order singularities}\end{center}

We begin with the fictitious Brownian motion model, in which the end-points of the paths are complex numbers.
The model consists of $2Rn$ nonintersecting Brownian paths starting at zero, with $n$ of them ending at each of the points $\pm n^{1/2}b_r\ (r=1,\ldots,R)$. The product in the integrand in (\ref{Hkl}) becomes
\be\prod_{r=1}^R\left({b_r^2-s^2/n\ov b_r^2-t^2/n}\right)^n,\label{Rprod}\ee
and we use the same contours as before. 

We shall first make the substitutions $\t_k\to 1/2+n^{-\dl}\t_k$ with $\dl$ to be determined. The first exponential in (\ref{Hkl}) becomes
\be\exp\Big\{-(1+4\,n^{-\dl}\t_k+O(n^{-2\dl}))\,t^2+(4+O(n^{-\dl}))\,xt+O(x^2)\Big\},\label{texp}\ee
and the second exponential becomes
\[\exp\Big\{(1+4\,n^{-\dl}\t_\l+O(n^{-2\dl}))\,s^2-(4+O(n^{-\dl}))\,ys+O(y^2)\Big\}.\]
If we set $a_r=1/b_r^2$ the product (\ref{Rprod}) is the exponential of
\[(\sum a_r)(t^2-s^2)+{n\inv\ov2}(\sum a_r^2)(t^4-s^4)+\cdots+{n^{-R}\ov R+1}(\sum a_r^{R+1})(t^{2R+2}-s^{2R+2})\]
\[+O(n^{-R-1}(|t|^{2R+4}+|s|^{2R+4})).\]

If $R>1$ we choose the $a_r$ such that
\[\sum a_r=1,\ \ \ \sum a_r^2=\cdots=\sum a_r^R=0.\]
The $a_r$ are the roots of the equation
\be a^R-a^{R-1}+{1\ov 2!}a^{R-2}-\cdots+{(-1)^R\ov R!}=0,\label{aeq}\ee
from which it follows that
\[\sum a_r^{R+1}={(-1)^{R+1}\ov R!}.\]
In general the $a_r$ will be complex, and so the same will be true of the end-points of our Brownian paths.

In the integrals defining $H_{ij}$ we make the substitutions
\[t\to n^{\dl/2}\,t,\ \ \ s\to n^{\dl/2}\,s,\]
where 
\[\dl={R\ov R+1},\]
and in the kernel
we make the scaling \[x\to n^{-\dl/2}\,x,\ \ \ y\to n^{-\dl/2}\,y.\]
This gives us the kernel $H_{n,ij}(x,y)$. As before $E_{ij}$ is unchanged, and the limiting form of $H_{n,ij}(x,y)$ is now
\[-{1\ov \pi^2}\int_C \int_{-i\iy}^{i\iy}e^{(-1)^Rs^{2R+2}/(R+1)!+4\,\t_j\,s^2-4ys+
(-1)^{R+1}t^{2R+2}/(R+1)!-4\,\t_i\,t^2+4xt}\,{ds\,dt\ov s-t}.\]

This is formal and it is not at first clear what the $C$ contour should be, although one might guess that it consists of four rays, one in each quadrant, on which $(-1)^{R+1}t^{2R+2}$ is negative and real. We shall see that this is so, and that the rays are the most vertical ones, those between 0 and $\pm\iy e^{\pm{R\ov2R+2}\pi i}$. The orientation of the rays is as in the case $R=1$. (The $s$ integration should cause no new problems.) 

After the variable changes the product of the two functions of $t$ in the integrand in (\ref{Hkl}) is of the form
\be {e^{-n^\dl t^2}\ov\prod (b_r^2-n^{\dl-1}t^2)^n}\,e^{(-4\,\t_k\,t^2+4xt+O(n^{-\dl}(|t|+|t|^2))}.\label{tintegrand}\ee
The main part of this is the quotient.

Upon the substitution $t\to n^{(\dl-1)/2}t$ the quotient becomes
\be{e^{-n t^2}\ov\prod (b_r^2-t^2)^n}.\label{tapprox}\ee
Suppose we want to do a steepest descent analysis of the integral of this
over a nearly vertical ray from 0 in the right half-plane. (This nearly vertical ray would be the part of the $t$ contour in (\ref{Hkl}) in the first quadrant.) No pole $\pm b_r$ is purely imaginary, as is clear from the equation the $a_r$ satisfy. So there are $R$ poles in the right half plane and $R$ in the left. There are $2R+2$ steepest descent curves emanating from the origin, half starting out in the right half-plane. These remain there since, as one can show,  the integrand is positive and increasing on the imaginary axis. We claim that there is at least one pole between any two of these curves. The reason is that otherwise the integrals over these two curves would be equal, and so have equal asymptotics. That means, after computing the asymptotics, that the integrals
\[\int_0^{\iy e^{ik\pi/(R+1)}} e^{(-1)^{R+1}\,t^{2R+2}}\,dt\]
would be the same for two different integers $k\in (-(R+1)/2,\,(R+1)/2)$. But the integrals are all different. 

Therefore there is a pole between any two of the curves.
Let $\G$ be the curve which starts out most steeply, in the direction ${\rm
arg}\;t={R\ov2R+2}\pi$. It follows from what we have just shown that there is no pole between $\G$ and the positive imaginary axis. This is what we wanted to show.

The curve we take for the $t$ integral in (\ref{Hkl}) is $\G_n=n^{(1-\dl)/2}\,\G$. The original contour for the $t$-integral in the representation of $H_{n,ij}$ can be deformed to this one. (We are speaking now, of course, of one quarter of the full contour.)

We can now take care of the annoying part of the argument establishing the claimed asymptotics. The curve $\G$ is asymptotic to the positive real axis at $+\iy$. Therefore for $A$ sufficiently large (\ref{tapprox}) is $O(e^{-n|t|^2/2})$ when $t\in\G,\ |t|>A$. Hence (\ref{tintegrand}) is $O(e^{-n^{\dl} |t|^2/4})$ when $t\in \G_n,\ |t|>n^{(1-\dl)/2} A$. It follows that its $L^2$ norm over this portion of $\G_n$ is exponentially small. When $t\in \G,\ |t|>\ep$ then (\ref{tapprox}) is $O(e^{-n\eta})$ for some $\eta>0$, and it follows that (\ref{tintegrand}) is $O(e^{-n\eta/2})$ when $t\in \G_n,\ |t|>n^{(1-\dl)/2} \ep$ and also $|t|<n^{(1-\dl)/2} A$. Therefore the norm of (\ref{tintegrand}) over this portion of $\G_n$ is also exponentially small. So we need consider only the portion of $\G_n$ on which $|t|<n^{(1-\dl)/2} \ep$, and for this we get the limit
\[\int _0^{\iy e^{2i\pi R/(R+1)}} e^{(-1)^{R+1}\,t^{2R+2}-4\,\t_k\,t^2+x\,t}\,dt\]
with appropriate uniformity, in the usual way.

Just as with the Pearcey kernel we can search for a system of PDEs associated with $\det\,(I-K\,\ch)$. Again we obtain two commutators, which when combined show that $K$ is an integrable kernel when $m=1$. In this case we define $\ph$ and $\ps$ by
\[\ph(x)=\left({1\ov \pi i}\int_Ce^{(-1)^{R+1}t^{2R+2}/(R+1)!-4\t_k\,t^2/2+4xt}dt\right),\]
\[\ps(y)=\left({1\ov \pi i}\int_{-i\iy}^{i\iy}
e^{(-1)^R s^{2R+2}/(R+1)!+4\t_k\,t^2/2-4yt}ds\right).\]
They satisfy the differential equations
\[c_R\,\ph^{(2R+1)}(x)-2\t\ph'(x)+4x\ph(x)=0,\ \ \ 
c_R\,\ps^{(2R+1)}(y)-2\ps'(y)\t-4y\ps(y)=0,\]
where
\[c_R=2{(-1)^{R+1}\ov 4^{2R+1}\,R!}.\]

The first commutator is
$[D,\,K]=-4\,\ph\tn\ps$, which also gives 
\be [D^n,\,K]=-4\sum_{k=0}^{n-1}(-1)^k\,\ph^{(n-k-1)}\tn\ps^{(k)}.\label{Dncom}\ee

The second commutator is
\[[c_R\,D^{2R+1}-2\,\t\,D+4M,\,K]=0.\]
In case $m=1$ (or for a general $m$ and a diagonal entry of $K$), combining this with  the commutator $[D^{2R+1},\,K]$ obtained from (\ref{Dncom}) and the differential equations for $\ph$ and $\ps$  we get an expression for $[M,\,K]$ in terms of derivatives of $\ph$ and $\ps$ up to order $2R$. This gives the analogue of the expression for $K(x,y)$ in footnote \ref{Krep1}.

For a system of PDEs in this case we would have many more unknowns, and the industrious reader could write them down. However there will remain the problem of showing that certain quantities involving $2R$th derivatives of the resolvent kernel $R$ (too many $R$s!) evaluated at endpoints of the $X_k$ are expressible in terms of the unknowns. For the case $R=1$ a miracle took place. Even to determine what miracle has to take place for general $R$ would be a nontrivial computational task.

\begin{center}{\bf Acknowledgments}\end{center}

This work was supported by the National Science Foundation under grants DMS-0304414 (first
author) and DMS-0243982 (second author).

\end{document}